\documentclass{amsart}
\title
{$N$-tuple sum analogues for Ramanujan-type congruences}
\usepackage{amssymb,amsmath,amsthm,epsfig,graphics,latexsym}
\usepackage{enumerate}
 \theoremstyle{definition}

  \theoremstyle{plain}
  \newtheorem{lemma}      {Lemma}
  
  \newtheorem{theorem}    {Theorem}
  
  \newtheorem{corollary}  {Corollary}
  
  \newtheorem*{conjecture} {Conjecture}

  \newcommand{\fr}{\frac}

\begin{document}
  \author{Mohamed El Bachraoui}
\email{melbachraoui@uaeu.ac.ae}
\keywords{$q$-analogue; super(congruence); cyclotomic polynomial; truncated sum}
\subjclass{11B65; 11F33; 33C20}
\begin{abstract}
In this paper we prove supercongruence relations for truncated $N$-tuple sums of basic hypergeometric series.
As an application, we give double, triple, and quadruple sum analogues of some Ramanujan-type supercongruences.
\end{abstract}
\date{\textit{\today}}
 \maketitle
\section{Introduction}
\noindent
Throughout the paper $n$ is a nonnegative integer, $p$ is a prime number,
and $x$ and $q$ are complex numbers such
that $|q|<1$. We will use the following notation from the theory of basic hypergeometric
series, see Gasper-Rahman~\cite{Gasper-Rahman}.
The $q$-shifted factorials are given by
\[
(x;q)_0= 1,\ \quad (x;q)_n = \prod_{j=0}^{n-1} (1-xq^j),\ \text{and\ }
(x;q)_{\infty} = \prod_{j=0}^{\infty} (1-xq^j).
\]
The $q$-integer is given for any nonnegative integer $n$ by
\[
[n]_q=[n] = \fr{1-q^n}{1-q} = 1+q+\cdots+q^{n-1}
\]
and it is easy to check that
$\lim_{q\to 1} [n]=n$.
Let $\mathbb{Z}$, $\mathbb{N}_0$, and $\mathbb{N}$ denote the sets of integers,
nonnegative integers, and positive integers respectively.
Let
$\mathbb{Z}[q]$ denote the set of polynomials in $q$ with integer coefficients and let
$\mathbb{Z}(q)$ denote the set of rational functions in $q$ with integer coefficients.
The $n$-th cyclotomic polynomial is the polynomial in $\mathbb{Z}[q]$ given by
\[
\Phi_n(q) = \prod_{\substack{j=1 \\ \gcd(j,n)=1}}^n(q- \zeta^j) ,
\]
where $\zeta = e^{2\pi i/n}$ is the $n$-th root of unity.
Given polynomials $A_1(q), A_2(q), P(q) \in \mathbb{Z}[q]$, the congruence
$\fr{A_1(q)}{A_2(q)} \equiv 0 \pmod{P(q)}$ means that $P(q)$ divides $A_1(q)$ and that
$\gcd\big(P(q),A_2(q) \big) =1$; and in general for rational functions $A(q), B(q) \in\mathbb{Z}(q)$,
the congruence $A(q)\equiv B(q) \pmod{P(q)}$ means that $A(q)-B(q)\equiv 0\pmod{P(q)}$.
Ramanujan-type supercongruences for truncated sums of basic hypergeometric series
have received much interest in recent years, see for
instance~\cite{Guillera-Zudilin, Guo-1, Guo-2, Guo-3, Guo-4, Guo-Wang, Guo-Zudilin-2, Guo-Zudilin-3, Li, Osburn-Zudilin, Song-Wang, Zudilin}.
Most of the publications dealing with supercongruence relations were motivated by Ramanujan's formulas for $\fr{1}{\pi}$ as
in Ramanujan~\cite{Ramanujan} and Berndt~\cite{Berndt} and the later work by Van~Hamme~\cite{VanHamme} where he
stated a variety of conjectures on the $p$-adic analogues of Ramanujan-type formulas for $1/\pi$.
Among Van Hamme's $p$-adic conjectures one finds
\[
\begin{split}
\text{(E2)\quad } & \sum_{k=0}^{p-1} (-1)^k \fr{\big(\fr{1}{3}\big)_{k}^3}{k! ^3} (6k+1) \equiv p \pmod{p^3}
\quad \text{for\ } p \equiv 1\pmod{3} \\
\text{(F2)\quad } & \sum_{k=0}^{p-1} (-1)^k \fr{\big(\fr{1}{4}\big)_{k}^3}{k! ^3} (8k+1) \equiv p \Big( \fr{-2}{p}\Big) \pmod{p^3}
\quad \text{for\ }  p \equiv 1\pmod{4}\\
\text{(G2)\quad } & \sum_{k=0}^{p-1}  \fr{\big(\fr{1}{4}\big)_{k}^4}{k! ^4} (8k+1)
\equiv p \fr{\big(\fr{1}{2}\big)_{(p-1)/4}}{\big(1\big)_{(p-1)/4}} \pmod{p^3}
\quad \text{for\ } p\equiv 1\pmod{4}
\end{split}
\]
where $(x)_n = x(x+1)\cdots (x+n-1)$ is the Pochhammer function and $\big(\fr{\cdot}{p}\big)$ is the Legendre symbol modulo $p$.
All three conjectures (E2), (F2), and (G2) have been settled by Swisher~\cite{Swisher}.
Guo and Zudilin~\cite{Guo-Zudilin-2} recently developed the \emph{$q$-microscope} approach and
proved a large number of supercongruences for truncated sums of basic hypergeometric series
along with a $q$-analogue for the following Ramanujan-type congruence
\begin{equation*}
\text{(GZ)\quad } \sum_{k=0}^{p-1} \fr{{4k \choose 2k}{2k\choose k}^2}{2^{8k} 3^{2k}} (8k+1) \equiv p \Big( \fr{-3}{p}\Big) \pmod{p^3} \quad \text{for\ } p>3.
\end{equation*}
Another example of this type of supercongruences which is important for our work is
\[
\text{(LW)\quad } \sum_{k=0}^{p-1} \fr{\big(\fr{1}{4}\big)_{k}^4}{k! ^4} (8k+1)
\equiv - p \fr{\big(\fr{1}{2}\big)_{(p-1)/4}}{\big(1\big)_{(p-1)/4}} \pmod{p^3}
\quad \text{for\ } p\equiv 1\pmod{4}
\]
which is recently obtained by Liu and Wang's~\cite[p. 3]{Liu-Wang} using the $q$-microscope method.
It is natural to raise the question about $N$-tuple sum analogues of the Ramanujan-type congruences.
In~\cite{Bachraoui} we dealt with the case $N=2$ as we established double sum versions of some Ramanujan-type congruences
including the following double sum analogue of~(GZ)
\[
\sum_{k=0}^{p-1} \sum_{j=0}^k \fr{{4j \choose 2j}{2j\choose j}^2
{4(k-j) \choose 2(k-j)}{2(k-j)\choose k-j}^2}{2^{8k} 3^{2k}} (8j+1)\big(8(k-j)+1 \big)
\equiv p^2 \pmod{p^3}.
\]
To achieve this, we combined the $q$-microscope approach with a result on double sums of series satisfying certain properties.
Recently, Li~\cite{Li} and Song and Wang~\cite{Song-Wang} gave other supercongruences for truncated sums of squares of $q$-series.
In this note we consider the general case where $N\geq 2$. As an application we deduce supercongruences for truncated $N$-tuple sums
of some (basic) hypergeometric series including double, triple, and quadruple sum analogues of
the relations (E2), (F2), (G2), (GZ), and (LW).
Proofs of all our main results rely on the $q$-microscope technique and the following lemma.
\begin{lemma}\label{lem:0}
Let $N\geq 2$ be a positive integer, let
$\{z_1(k)\}_{k=0}^{\infty},\ldots, \{z_N(k)\}_{k=0}^{\infty}$ be sequences of complex numbers
and let
\[
t(k) = \sum_{i_1 =0}^{k} \sum_{i_2 =0}^{k-i_1}\cdots \sum_{i_{N-1}=0}^{k-i_1-\cdots-i_{N-1}}
z_{1}(i_1)\cdots z_{N-1}(i_{N-1}) z_N (k-i_1-\cdots-i_{N-1}).
\]
(a)\ If there is $d\in\mathbb{N}$ such that for all $j=1,\ldots,N$
\begin{equation}\label{Eq:key-1}
z_j(k)= 0\ \text{for all $k\in\mathbb{N}_0$\ such that\ } \fr{d-1}{N} <k <d,
\end{equation}
then
\[
\sum_{i=1}^N \sum_{k=0}^{d-1} z_i (k)
= \sum_{k=0}^{d-1} t(k).
\]
(b)\ If moreover,
\begin{equation}\label{Eq:key-2}
\begin{split}
z_j (ld+k) &= z_j(ld) z_j (k)\ \text{for all $k,l\in\mathbb{N}_0$ with $0\leq k< d$},
\end{split}
\end{equation}
then for any $k, l\in\mathbb{N}_0$ such that $0\leq k <d$, there holds
\[
t(ld+k) = t(k) \sum_{i_1 =0}^{l} \sum_{i_2 =0}^{l-i_1}\cdots \sum_{i_{N}=0}^{l-i_1-\cdots-i_{N-1}}
z_{1}(i_1d)\cdots z_{N-1}(i_{N-1}d) z_N \big( (l-i_1-\cdots-i_{N-1})d \big).
\]
(c)\ If for all $j=1,\ldots,N$, the sequence $\{z_j(k)\}_{k=0}^{\infty}$ satisfies both~\eqref{Eq:key-1}-\eqref{Eq:key-2} for some $d$
and if $n$ is a positive multiple of $d$, then
\[
\sum_{k=0}^{n-1} t(k)=
\Big( \sum_{l=0}^{\fr{n}{d}-1} \sum_{i_1 =0}^{l} \sum_{i_2 =0}^{l-i_1}\cdots \sum_{i_{N}=0}^{l-i_1-\cdots-i_{N-1}}
z_{1}(i_1d)\cdots z_{N-1}(i_{N-1}d) z_N \big( (l-i_1-\cdots-i_{N-1})d \big) \Big)
\]
\[
\times \sum_{k=0}^{d-1} t(k).
\]
\end{lemma}
The paper is organized as follows. In Section~\ref{sec:results} we state our main theorems along with
their consequences. Proofs for these theorems are given in Sections~\ref{sec:appl-2-proof}-\ref{sec:appl-5-proof}.
Section~\ref{sec:lem0-proof} is devoted to the proof of Lemma~\ref{lem:0}. In Section~\ref{sec:comments} we close with some comments and suggestions for further research.
\section{Statement of results}\label{sec:results}
\begin{theorem}\label{thm:appl-1}
Let $s, N\in\mathbb{N}$ such that $2\leq N \leq s$, let
\[
\alpha_q(k) = (-1)^k q^{s{k+1\choose 2}-k} [2sk+1] \fr{(xq;q^s)_k (q/x;q^s)_k (q;q^s)_k}
{(xq^s;q^s)_k (q^s/x;q^s)_k (q^s;q^s)_k},
\]
and let
\[
t_q(k) = \sum_{i_1=0}^k \sum_{i_2=0}^{k-i_1}\cdots \sum_{i_N =0}^{k-i_1-\cdots i_{N-1}}
\alpha_q(i_1) \cdots \alpha_q(i_{N-1}) \alpha_q(k-i_1- \cdots-i_{N-1}).
\]
Then for any positive integer $n$ whose prime factors are all congruent to $1$ modulo $s$,
there holds
\[
\sum_{k=0}^{n-1} t_q(k)
\equiv q^{N(n-1)(n-s+1)/(2s)}[n]^N (-1)^{N(n-1)/s}
\pmod{[n](1- a q^n)(a-q^n)}.
\]
\end{theorem}
In our applications of Theorem~\ref{thm:appl-1} we restrict ourselves to $s=3$ and $s=4$
as they are related to Van Hamme's relations (E2) and (F2).
So, letting in Theorem~\ref{thm:appl-1}, $s=3$, $x=1$ and $q\to 1$, we obtain the
following double and triple sum analogues (E2).
\begin{corollary}\label{cor:appl-1-1}
If $p\equiv 1 \pmod{3}$, then
\[
\sum_{k=0}^{p-1}\sum_{j=0}^k (-1)^k \fr{\big(\fr{1}{3}\big)_{j}^3}{j! ^3}
\fr{\big(\fr{1}{3}\big)_{k-j}^3}{(k-j)! ^3}
(6j+1)\big(6(k-j)+1\big) \equiv p^2  \pmod{p^3}
\]
and
\[
\sum_{k=0}^{p-1}\sum_{j=0}^k \sum_{i=0}^{k-j} (-1)^k \fr{\big(\fr{1}{3}\big)_{j}^3}{j! ^3}
\fr{\big(\fr{1}{3}\big)_{i}^3}{i! ^3} \fr{\big(\fr{1}{3}\big)_{k-j-i}^3}{(k-j-i)! ^3}
(6i+1)(6j+1)\big(6(k-j-i)+1 \big) \equiv 0  \pmod{p^3}.
\]
\end{corollary}
Letting in Theorem~\ref{thm:appl-1}, $s=4$, $x=1$ and $q\to 1$, we deduce the following
double, triple, and quadruple sum analogues of (F2).
\begin{corollary}\label{cor:appl-1-2}
If $p\equiv 1 \pmod{4}$, then
\[
\sum_{k=0}^{p-1}\sum_{j=0}^k  (-1)^k \fr{\big(\fr{1}{4}\big)_{j}^3}{j! ^3}
 \fr{\big(\fr{1}{4}\big)_{k-j}^3}{(k-j)! ^3}
(8j+1) \big(8(k-j)+1 \big) \equiv p^2  \pmod{p^3}
\]
and
\[
\sum_{k=0}^{p-1}\sum_{j=0}^k \sum_{i=0}^{k-j} (-1)^k \fr{\big(\fr{1}{4}\big)_{j}^3}{j! ^3}
\fr{\big(\fr{1}{4}\big)_{i}^3}{i! ^3} \fr{\big(\fr{1}{4}\big)_{k-j-i}^3}{(k-j-i)! ^3}
(8i+1)(8j+1) \big(8(k-j-i)+1 \big)
\]
\[
\equiv \sum_{k=0}^{p-1}\sum_{j=0}^k \sum_{i=0}^{k-j}\sum_{h=0}^{k-j-i} (-1)^k
\fr{\big(\fr{1}{4}\big)_{j}^3}{j! ^3} \fr{\big(\fr{1}{4}\big)_{i}^3}{i! ^3}
\fr{\big(\fr{1}{4}\big)_{h}^3}{h! ^3} \fr{\big(\fr{1}{4}\big)_{k-j-i-h}^3}{(k-j-i-h)! ^3}
(8i+1)(8j+1)(8h+1)
\]
\[
\qquad \times \big(8(k-j-i-h)+1\big) \equiv 0  \pmod{p^3}.
\]
\end{corollary}
Theorem~\ref{thm:appl-1} can be generalized as follows.
\begin{theorem}\label{thm:appl-2}
Let $(r,s)\in\mathbb{Z}\times\mathbb{N}$ such that $\gcd(r,s)=1$, let
\[
\alpha_q(k) = (-1)^k q^{s{k+1\choose 2}-rk} [2sk+1] \fr{(xq^r;q^s)_k (q^r/x;q^s)_k (q^r;q^s)_k}
{(xq^s;q^s)_k (q^s/x;q^s)_k (q^s;q^s)_k},
\]
and let for a positive integer $N\geq 2$
\[
t_q(k) = \sum_{i_1=0}^k \sum_{i_2=0}^{k-i_1}\cdots \sum_{i_N =0}^{k-i_1-\cdots i_{N-1}}
\alpha_q(i_1) \cdots \alpha_q(i_{N-1}) \alpha_q(k-i_1- \cdots-i_{N-1}).
\]
Then for any prime  $p\equiv r\pmod{s}$ such that $r\leq p$ and $(p-r)/s \leq (p-1)/N$,
there holds
\begin{equation}\label{r-s}
\sum_{k=0}^{p-1} t_q(k)
\equiv q^{N(p-r)(p-s+r)/(2s)}[p]^N (-1)^{N(p-r)/s}
\pmod{[p](1- a q^p)(a-q^p)}.
\end{equation}
\end{theorem}
Theorem~\ref{thm:appl-2} obviously has a wide range of applications, but we restrict ourselves to the following
two cases.
Letting in Theorem~\ref{thm:appl-2}, $(r,s)=(2,3)$, $x=1$ and $q\to 1$, we get.
\begin{corollary}
If $p\equiv 2 \pmod{3}$, then
\[
\sum_{k=0}^{p-1}\sum_{j=0}^k (-1)^k \fr{\big(\fr{2}{3}\big)_{j}^3}{j! ^3}
 \fr{\big(\fr{2}{3}\big)_{k-j}^3}{(k-j)! ^3}
(6j+1)\big(6(k-j)+1\big) \equiv p^2  \pmod{p^3}
\]
and
\[
\sum_{k=0}^{p-1}\sum_{j=0}^k \sum_{i=0}^{k-j} (-1)^k \fr{\big(\fr{2}{3}\big)_{j}^3}{j! ^3}
\fr{\big(\fr{2}{3}\big)_{i}^3}{i! ^3} \fr{\big(\fr{2}{3}\big)_{k-j-i}^3}{(k-j-i)! ^3}
(6i+1)(6j+1)\big(6(k-j-i)+1\big) \equiv 0  \pmod{p^3}.
\]
\end{corollary}
Upon letting in Theorem~\ref{thm:appl-2}, $(r,s)=(3,4)$, $x=1$ and $q\to 1$, we derive.
\begin{corollary}
If $p\equiv 3 \pmod{4}$, then
\[
\sum_{k=0}^{p-1}\sum_{j=0}^k (-1)^k \fr{\big(\fr{3}{4}\big)_{j}^3}{j! ^3}
 \fr{\big(\fr{3}{4}\big)_{k-j}^3}{(k-j)! ^3}
(8j+1) \big(8(k-j)+1 \big) \equiv p^2  \pmod{p^3}
\]
and
\[
\sum_{k=0}^{p-1}\sum_{j=0}^k \sum_{i=0}^{k-j} (-1)^k \fr{\big(\fr{3}{4}\big)_{j}^3}{j! ^3}
\fr{\big(\fr{3}{4}\big)_{i}^3}{i! ^3} \fr{\big(\fr{3}{4}\big)_{k-j-i}^3}{(k-j-i)! ^3}
(8i+1)(8j+1) \big(8(k-j-i)+1 \big)
\]
\[
\equiv \sum_{k=0}^{p-1}\sum_{j=0}^k \sum_{i=0}^{k-j}\sum_{h=0}^{k-j-i} (-1)^k
\fr{\big(\fr{3}{4}\big)_{j}^3}{j! ^3} \fr{\big(\fr{3}{4}\big)_{i}^3}{i! ^3}
\fr{\big(\fr{3}{4}\big)_{h}^3}{h! ^3} \fr{\big(\fr{3}{4}\big)_{k-j-i-h}^3}{(k-j-i-h)! ^3}
(8i+1)(8j+1)(8h+1)
\]
\[
\qquad \times \big(8(k-j-i-h)+1\big) \equiv 0  \pmod{p^3}.
\]
\end{corollary}
\begin{theorem}\label{thm:appl-3}
Let
\[
\alpha_q(k) = q^{2k} [8k+1] \fr{(q;q^4)_k^2 (xq;q^4)_k (q/x;q^4)_k}
{(q^4;q^4)_k^2 (xq^4;q^4)_k (q^4/x;q^4)_k},
\]
\[
s_q(k) = \sum_{j=0}^k \alpha_q(j)  \alpha_q(k-j), \quad
t_q(k) = \sum_{j=0}^k\sum_{i=0}^{k-j} \alpha_q(j) \alpha_q(i) \alpha_q(k-j-i),
\]
and
\[
u_q(k) = \sum_{j=0}^k\sum_{i=0}^{k-j} \sum_{h=0}^{k-j-i} \alpha_q(j) \alpha_q(i) \alpha_q(h)
\alpha_q(k-j-i-h).
\]
Then for any positive integer $n$ whose prime factors are all congruent to $1$ modulo $4$,
there holds
\[
\sum_{k=0}^{n-1} s_q(k)
\equiv [n]^2 \fr{(q^2;q^4)_{(n-1)/4}^2}{(q^4;q^4)_{(n-1)/4}^2}  q^{(1-n)/2}
\pmod{[n](1- x q^n)(x-q^n)},
\]
\[
\sum_{k=0}^{n-1} t_q(k)
\equiv [n]^3 \fr{(q^2;q^4)_{(n-1)/4}^3}{(q^4;q^4)_{(n-1)/4}^3}  q^{3(1-n)/4}
\pmod{[n](1- x q^n)(x-q^n)},
\]
and
\[
\sum_{k=0}^{n-1} u_q(k)
\equiv [n]^4 \fr{(q^2;q^4)_{(n-1)/4}^4}{(q^4;q^4)_{(n-1)/4}^4}  q^{(1-n)}
\pmod{[n](1- x q^n)(x-q^n)}.
\]
\end{theorem}
Upon letting in Theorem~\ref{thm:appl-3}, $x=1$ and $q\to 1$, we get the following double, triple, and quadruple sum versions of Van Hamme's congruence (G2).
\begin{corollary}\label{cor:appl-3-1}
If $p\equiv 1 \pmod{4}$ is a prime number, then we have
\[
\sum_{k=0}^{p-1}\sum_{j=0}^k \fr{\big(\fr{1}{4}\big)_{j}^4}{j! ^4}
 \fr{\big(\fr{1}{4}\big)_{k-j}^4}{(k-j)! ^4}
(8j+1)(8k-8j+1) \equiv p^2 \fr{\big(\fr{1}{2}\big)_{(p-1)/4}^2}{\big(1\big)_{(p-1)/4}^2} \pmod{p^3}
\]
and
\[
\sum_{k=0}^{p-1}\sum_{j=0}^k \sum_{i=0}^{k-j} \fr{\big(\fr{1}{4}\big)_{j}^4}{j! ^4}
\fr{\big(\fr{1}{4}\big)_{i}^4}{i! ^4} \fr{\big(\fr{1}{4}\big)_{k-j-i}^4}{(k-j-i)! ^4}
(8i+1)(8j+1)(8k-8j-8i+1)
\]
\[
\equiv
\sum_{k=0}^{p-1}\sum_{j=0}^k \sum_{i=0}^{k-j} \sum_{h=0}^{k-j-i} \fr{\big(\fr{1}{4}\big)_{j}^4}{j! ^4}
\fr{\big(\fr{1}{4}\big)_{i}^4}{i! ^4} \fr{\big(\fr{1}{4}\big)_{h}^4}{h! ^4}
\fr{\big(\fr{1}{4}\big)_{k-j-i-h}^4}{(k-j-i-h)! ^4} (8i+1)(8j+1)(8h+1)
\]
\[
\times  \big( 8(k-j-i-h)+1\big) \equiv 0  \pmod{p^3}.
\]
\end{corollary}
\begin{theorem}\label{thm:appl-4}
Let
\[
\alpha_q(k) = q^{-4k} [8k+1]^2 [8k+1]_{q^2} \fr{(q^2;q^8)_k^2 (xq^2;q^8)_k (q^2/x;q^8)_k}
{(q^8;q^8)_k^2 (xq^8;q^8)_k (q^8/x;q^8)_k},
\]
\[
s_q(k) = \sum_{j=0}^k \alpha_q(j)  \alpha_q(k-j), \quad
t_q(k) = \sum_{j=0}^k\sum_{i=0}^{k-j} \alpha_q(j) \alpha_q(i) \alpha_q(k-j-i),
\]
and
\[
u_q(k) = \sum_{j=0}^k\sum_{i=0}^{k-j} \sum_{h=0}^{k-j-i} \alpha_q(j) \alpha_q(i) \alpha_q(h)
\alpha_q(k-j-i-h).
\]
Then for any positive integer $n$ whose prime factors are all congruent to $1$ modulo $4$, modulo $[n]_{q^2}(1- x q^{2n})(x-q^{2n})$
there holds
\[
\sum_{k=0}^{n-1} s_q(k)
\equiv [n]_{q^2}^2 \fr{(q^4;q^8)_{(n-1)/4}^2}{(q^8;q^8)_{(n-1)/4}^2}  q^{(1-n)}
\Big( 1 - \fr{(1-xq^2)(1-q^2/x)}{(1-q)^2 (1+q^2)} \Big)^2,
\]
\[
\sum_{k=0}^{n-1} t_q(k)
\equiv [n]_{q^2}^3 \fr{(q^4;q^8)_{(n-1)/4}^3}{(q^8;q^8)_{(n-1)/4}^3}  q^{3(1-n)/2}
\Big( 1 - \fr{(1-xq^2)(1-q^2/x)}{(1-q)^2 (1+q^2)} \Big)^3,
\]
and
\[
\sum_{k=0}^{p-1} u_q(k)
\equiv [n]_{q^2}^4 \fr{(q^4;q^8)_{(n-1)/4}^4}{(q^8;q^8)_{(n-1)/4}^4}  q^{2(1-n)}
\Big( 1 - \fr{(1-xq^2)(1-q^2/x)}{(1-q)^2 (1+q^2)} \Big)^4.
\]
\end{theorem}
Letting in Theorem~\ref{thm:appl-4}, $x=1$ and $q\to 1$, we derive the following double, triple, and quadruple versions of~(LW).
\begin{corollary}\label{cor:appl-4-1}
If $p\equiv 1 \pmod{4}$ is a prime number, then we have
\[
\sum_{k=0}^{p-1}\sum_{j=0}^k \fr{\big(\fr{1}{4}\big)_{j}^4}{j! ^4}
 \fr{\big(\fr{1}{4}\big)_{k-j}^4}{(k-j)! ^4}
(8j+1)^3(8k-8j+1)^3 \equiv p^2 \fr{\big(\fr{1}{2}\big)_{(p-1)/4}^2}{\big(1\big)_{(p-1)/4}^2} \pmod{p^3}
\]
and
\[
\sum_{k=0}^{p-1}\sum_{j=0}^k \sum_{i=0}^{k-j} \fr{\big(\fr{1}{4}\big)_{j}^4}{j! ^4}
\fr{\big(\fr{1}{4}\big)_{i}^4}{i! ^4} \fr{\big(\fr{1}{4}\big)_{k-j-i}^4}{(k-j-i)! ^4}
(8i+1)^3(8j+1)^3(8k-8j-8i+1)^3
\]
\[
\equiv
\sum_{k=0}^{p-1}\sum_{j=0}^k \sum_{i=0}^{k-j}  \sum_{h=0}^{k-j-i}\fr{\big(\fr{1}{4}\big)_{j}^4}{j! ^4}
\fr{\big(\fr{1}{4}\big)_{i}^4}{i! ^4} \fr{\big(\fr{1}{4}\big)_{h}^4}{h! ^4}
\fr{\big(\fr{1}{4}\big)_{k-j-i-h}^4}{(k-j-i-h)! ^4}
(8i+1)^3(8j+1)^3 (8h+1)^3
\]
\[
\qquad\qquad\qquad \times \big(8(k-j-i-h)+1\big)^3 \equiv 0  \pmod{p^3}.
\]
\end{corollary}
Our last application is related to the congruence~(GZ).
\begin{theorem}\label{thm:appl-5}
Let
\[
\alpha_q(k)= \fr{(xq;q^2)_k (q/x;q^2)_k (q;q^2)_{2k}}{(xq^6;q^6)_k (q^6/x;q^6)_k (q^2;q^2)_{2k}},
[8k+1] q^{2k^2},
\]
let
\[
s_q(k) = \sum_{j=0}^k \alpha_q(j)  \alpha_q(k-j), \quad
t_q(k) = \sum_{j=0}^k\sum_{i=0}^{k-j} \alpha_q(j) \alpha_q(i) \alpha_q(k-j-i),
\]
and
\[
u_q(k) = \sum_{j=0}^k\sum_{i=0}^{k-j} \sum_{h=0}^{k-j-i}\alpha_q(j) \alpha_q(i) \alpha_q(i)
\alpha_q(k-j-i-h).
\]
Then for any positive integer $n$ which is coprime to $6$ there holds
\[
\sum_{k=0}^{n-1} s_q(k) \equiv q^{1-n}[n]^2 \pmod{[n](1-xq^n)(x-q^n)}
\]
and
\[
\sum_{k=0}^{n-1} t_q(k) \equiv \sum_{k=0}^{n-1} u_q(k) \equiv 0 \pmod{[n]}.
\]
\end{theorem}
\noindent
Note that the supercongruence in Theorem~\ref{thm:appl-5} for the double sum $\sum_{k=0}^{n-1} s_q(k)$ already appeared in~\cite[Theorem 3]{Bachraoui} and we include it here
for completeness.
Upon letting $x = 1$ and then $q\to 1$ in Theorem~\ref{thm:appl-5} we obtain the following
triple and quadruple sum versions of~(GZ).
\begin{corollary}\label{cor:main-5}
If $p>3$ be a prime number, then we have
\[
\sum_{k=0}^{p-1} \sum_{j=0}^k \sum_{i=0}^{k-j}
\fr{ {2j\choose j}^2{4j\choose 2j} {2i\choose i}^2{4i\choose 2i}
{2(k-j-i)\choose k-j-i}^2{4(k-j-i)\choose 2(k-j-i)}}{2^{8k} 3^{2k}}
(8i+1)(8j+1) \big(8(k-j-i)+1 \big)
\]
\[
\equiv \sum_{k=0}^{p-1} \sum_{j=0}^k \sum_{i=0}^{k-j} \sum_{h=0}^{k-j-i}
\fr{ {2j\choose j}^2{4j\choose 2j} {2i\choose i}^2{4i\choose 2i} {2h\choose h}^2{4h\choose 2h}
{2(k-j-i-h)\choose k-j-i}^2 {4(k-j-i-h)\choose 2(k-j-i-h)}}{2^{8k} 3^{2k}}
\]
\[
\times (8i+1)(8j+1)(8h+1)\big(8(k-j-i-h)+1 \big) \equiv 0  \pmod{p}.
\]
\end{corollary}
\section{Proof of Theorem~\ref{thm:appl-2}}\label{sec:appl-2-proof}
\noindent
Let $p$ be a prime number such that $p\equiv r\pmod{s}$ and let $\zeta$ be a $p$-th root of unity.
From~\cite[p. 349]{Guo-Zudilin-2} we have
\begin{equation}\label{Eq:GZ-1-1}
\sum_{k=0}^{p-1} \alpha_{\zeta}(k) =0.
\end{equation}
Furthermore, note that from the assumptions $p\equiv r\pmod{s}$, $r\leq p$, and $(p-r)/s \leq (p-1)/N$ it follows that
\[
\alpha_{\zeta}(k)=0\ \text{for\ } \fr{p-1}{N}\leq k \leq p-1,
\]
implying that the sequence
$\{\alpha_{\zeta}(k)\}_{k=0}^{\infty}$ satisfies condition~\eqref{Eq:key-1}.
Then by an appeal to Lemma~\ref{lem:0}(a) with
$z_j (k)= \alpha_{\zeta}(k)$ for $j=1,\ldots,N$, we find
\[
\sum_{k=0}^{p-1} t_{\zeta}(k)
= \Big(\sum_{k=0}^{p-1} \alpha_{\zeta}(k) \Big) ^N
\]
and so by~\eqref{Eq:GZ-1-1},
\begin{equation}\label{help1-thm-1}
\sum_{k=0}^{p-1} t_{q}(k) \equiv 0 \equiv q^{N(p-r)(p-s+r)/(2s)}[p]^N (-1)^{N(p-r)/s}\pmod{[p]}.
\end{equation}
Furthermore, from~\cite[p. 349]{Guo-Zudilin-2} by letting $x=q^{-p}$ or $x=q^p$ in $\alpha_q(k)$, we
have
\begin{equation}\label{help2-thm-1}
\sum_{k=0}^{p-1} \alpha_q(k)
= q^{(p-r)(p-s+r)/(2s)} [p] (-1)^{(p-r)/s}.
\end{equation}
Noting that with these choices of $x$, the sequence $\{\alpha_q(k)\}_{k=0}^{\infty}$ fulfills
property~\eqref{Eq:key-1}, we deduce from Lemma~\ref{lem:0}(a) that
\begin{equation}\label{help3-thm-1}
\sum_{k=0}^{p-1} t_q(k)  = \Big( \sum_{k=0}^{p-1} \alpha_q(k) \Big)^N.
\end{equation}
As the polynomials $[p]$, $(1-x q^p)$, and $(x-q^p)$ are relatively prime,
we achieve the desired congruence by combining~\eqref{help1-thm-1}-\eqref{help3-thm-1}.
\section{Proof of Theorem~\ref{thm:appl-1}}\label{sec:appl-1-proof}
Let $d\mid n$ and let $\zeta$ be a $d$-th primitive root of unity.
From the assumptions that $d\equiv 1\pmod{s}$ and that $N\leq s$ we get
\[
\alpha_{\zeta}(k)= 0\ \text{for\ } \fr{d-1}{N} \leq k \leq d-1
\]
and by ~\cite{Guo-Zudilin-2} that
\[
\sum_{k=0}^{d-1} \alpha_{\zeta}(k) =0
\]
and so, by Lemma~\ref{lem:0}(a)
\[
\sum_{k=0}^{d-1} t_{\zeta}(k)  = \Big( \sum_{k=0}^{d-1} \alpha_{\zeta}(k) \Big)^N =0.
\]
Then the forgoing identity and Lemma~\ref{lem:0}(c) give
\[
\sum_{k=0}^{n-1} t_{\zeta}(k)
=
\sum_{k=0}^{d-1} t_{\zeta}(k) \Big(\sum_{l=0}^{n/d-1}
\sum_{m=0}^l \sum_{i=0}^{l-m} \alpha_{\zeta}(md) \alpha_{\zeta}(id)\alpha_{\zeta}\big( (l-m-i)d\big)\Big)
= 0,
\]
yielding that
\begin{equation}\label{help1-thm-2b}
\sum_{k=0}^{n-1} t_q(k) \equiv 0 \equiv q^{N(n-1)(n-s+1)/(2s)}[n]^N (-1)^{N(n-1)/s}\pmod{[n]}.
\end{equation}
In addition, from~\cite{Guo-Zudilin-2} by letting $x=q^{-n}$ or $x=q^n$ in $\alpha_q(k)$, we
have
\begin{equation}\label{help2-thm-2b}
\sum_{k=0}^{n-1} \alpha_q(k)
= q^{(n-1)(n-s+1)/(2s)} [n] (-1)^{(n-1)/s}
\end{equation}
and it easily verified that the sequence $\{\alpha_q(k)\}_{k=0}^{\infty}$ fulfills
property~\eqref{Eq:key-1} as $x=q^{-n}$ or $x=q^{n}$. Then we deduce by Lemma~\ref{lem:0}(a) that
\begin{equation}\label{help3-thm-2b}
\sum_{k=0}^{n-1} t_q(k)  = \Big( \sum_{k=0}^{n-1} \alpha_q(k) \Big)^N.
\end{equation}
Now use ~\eqref{help1-thm-2b}-\eqref{help3-thm-2b} and the fact the polynomials $[n]$, $(1-x q^n)$, and
$(x-q^n)$ are relatively prime to derive the desired congruence.
\section{Proof of Theorem~\ref{thm:appl-3} and Theorem~\ref{thm:appl-4}}\label{sec:appl-3-4-proof}
\emph{Proof of Theorem~\ref{thm:appl-3}.}
Let $n$ be a positive integer whose prime factors are all congruent to $1$ modulo $4$,
Let $d\mid n$, and let $\zeta$ be a $d$-th root of unity.  Clearly $n \equiv d \equiv 1 \pmod{4}$. Liu-Wang~\cite[p. 5]{Liu-Wang}
proved that
\begin{equation}\label{appl-3 help-1}
\sum_{k=0}^{d-1} \alpha_{\zeta}(k) = 0
\end{equation}
and for $x=q^{n}$ or $x=q^{-n}$,
\begin{equation}\label{appl-3 help-2}
\sum_{k=0}^{n-1} \alpha_{q}(k) = [n] \fr{(q^2;q^4)_{(n-1)/4}}{(q^4;q^4)_{(n-1)/4}}  q^{(1-n)/4}.
\end{equation}
It easy to see that the sequence $\{\alpha_{\zeta}(k) \}_{k=0}^\infty$ satisfies both
conditions~\eqref{Eq:key-1} and~\eqref{Eq:key-2} for $N=1,2,3$.
Then by Lemma~\ref{lem:0}(a) and~\eqref{appl-3 help-1}, we find
\[
\sum_{k=0}^{d-1} s_{\zeta}(k)
= \Big(\sum_{k=0}^{d-1} \alpha_{\zeta}(k) \Big) ^2
\sum_{k=0}^{d-1} t_{\zeta}(k)
= \Big(\sum_{k=0}^{d-1} \alpha_{\zeta}(k) \Big) ^3
= \sum_{k=0}^{d-1} u_{\zeta}(k)
= \Big(\sum_{k=0}^{d-1} \alpha_{\zeta}(k) \Big) ^4 =0.
\]
Thus with the help of Lemma~\ref{lem:0}(c)
\[
\sum_{k=0}^{n-1} s_{\zeta}(k) \sum_{k=0}^{n-1} t_{\zeta}(k) = \sum_{k=0}^{n-1} u_{\zeta}(k)
= 0,
\]
implying that
\[
\sum_{k=0}^{n-1} s_{q}(k) \equiv \sum_{k=0}^{n-1} t_{q}(k)\equiv \sum_{k=0}^{n-1} u_{q}(k) \equiv 0 \pmod{[n]}.
\]
Besides, noticing that for $x=q^n$ or $x=q^{-n}$, the sequence $\{\alpha_{q}(k) \}_{k=0}^\infty$ meets the
property~\eqref{Eq:key-1}, we deduce from Lemma~\ref{lem:0}(a) that
\[
\sum_{k=0}^{n-1} s_{q}(k) = \Big(\sum_{k=0}^{n-1} \alpha_{q}(k) \Big)^2,\
\sum_{k=0}^{n-1} t_{q}(k) = \Big(\sum_{k=0}^{n-1} \alpha_{q}(k) \Big)^3,\ \text{and\ }
\sum_{k=0}^{n-1} u_{q}(k) = \Big(\sum_{k=0}^{n-1} \alpha_{q}(k) \Big)^4.
\]
This together with~\eqref{appl-3 help-2} yield that modulo $(1-xq^n)(x-q^n)$,
\[
\begin{split}
\sum_{k=0}^{n-1} s_{q}(k) &\equiv [n]^2 \fr{(q^2;q^4)_{(n-1)/4}^2}{(q^4;q^4)_{(n-1)/4}^2}  q^{(1-n)/2} \\
\sum_{k=0}^{n-1} t_{q}(k) &\equiv [n]^3 \fr{(q^2;q^4)_{(n-1)/4}^3}{(q^4;q^4)_{(n-1)/4}^3}  q^{3(1-n)/4} \\
\sum_{k=0}^{n-1} u_{q}(k) &\equiv [n]^4 \fr{(q^2;q^4)_{(n-1)/4}^4}{(q^4;q^4)_{(n-1)/4}^4}  q^{1-n}.
\end{split}
\]
Now combine the above congruences to achieve the desired result.

\emph{Proof of Theorem~\ref{thm:appl-4}.} Follows similarly by a combination of Lemma~\ref{lem:0} with
the proof of~\cite[Theorem 4]{Liu-Wang}.
\section{Proof of Theorem~\ref{thm:appl-5}}\label{sec:appl-5-proof}
The result is clear for $n=1$. Let $n>1$ with $\gcd(n,6)=1$, let $d\mid n$, and let $\zeta$ be a $d$-th
primitive root of unity.
From~\cite[p. 339]{Guo-Zudilin-2} we have
\begin{equation}\label{Eq:GZ-0}
\sum_{k=0}^{d-1} \alpha_{\zeta}(k) =0.
\end{equation}
As $(q;q^2)_{2k} = (q;q^4)_k (q^3;q^4)_k$
appears in the numerator of $\alpha_q(k)$ we see that
the sequence $\{\alpha_{\zeta}(k)\}_{k=0}^{\infty}$ satisfies
condition~\eqref{Eq:key-1} for $N\in\{3,4\}$.
Then by~\eqref{Eq:GZ-0} and Lemma~\ref{lem:0}(a) applied to $z_j(k)=\alpha_q(k)$ for $j=1,2,3,4$, we obtain
\[
\sum_{k=0}^{d-1} t_{\zeta}(k)
= \Big(\sum_{k=0}^{d-1} \alpha_{\zeta}(k) \Big) ^3 =0
\ \text{and\ }
\sum_{k=0}^{d-1} u_{\zeta}(k)
= \Big(\sum_{k=0}^{d-1} \alpha_{\zeta}(k) \Big) ^4 =0.
\]
Moreover, it is easy to check by the basic properties of roots of unity that $\{\alpha_{\zeta}(k)\}_{k=0}^{\infty}$ fulfills the
property~\eqref{Eq:key-2} as well and that $\alpha_{\zeta}(0)=1$.
It follows by virtue of Lemma~\ref{lem:0}(c) with $N=3$ that
\[
\sum_{k=0}^{n-1} t_{\zeta}(k) = \sum_{k=0}^{d-1} t_{\zeta}(k) \Big(\sum_{l=0}^{n/d-1}
\sum_{m=0}^l \sum_{i=0}^{l-m} \alpha_{\zeta}(md) \alpha_{\zeta}(id)\alpha_{\zeta}\big( (l-m-i)d\big)\Big)
=0,
\]
implying that $\sum_{k=0}^{n-1} t_q(k)$ is divisible by the cyclotomic polynomial $\Phi_d(q)$ for any divisor
$d>1$ of $n$. Hence
\begin{equation}\label{help0-app-1}
\sum_{k=0}^{n-1} t_q(k) \equiv 0  \pmod{[n]}.
\end{equation}
We similarly show by an appeal to Lemma~\ref{lem:0}(b, c) with $N=4$ that
$\sum_{k=0}^{n-1} u_q(k) \equiv 0  \pmod{[n]}$. This completes the proof.
%
\section{Proof of Lemma~\ref{lem:0}}\label{sec:lem0-proof}
\noindent
Part (a) is a direct consequence of~\eqref{Eq:key-1}. Part (c) follows from part~(b) and the relation
\[
\sum_{k=0}^{n-1}t(k) = \sum_{l=0}^{\fr{n}{d}-1} \sum_{k=0}^{d-1} t(ld+k).
\]
As for part~(b), we use induction on $N$.
Let $N=2$ and for simplicity of notation, let $a(k)= z_1(k)$, $b(k)=z_2(k)$,
and $c(k)=\sum_{j=0}^k a(j)b(k-j)$.
Then we have
\begin{equation}\label{Eq:lem-2}
c(ld+k)
= \sum_{j=0}^{ld+k} a(j) b(ld+k-j)
= \sum_{j=0}^{d-1} a(j) b(ld+k-j)
\end{equation}
\[
+ \sum_{j=d}^{2d-1} a(j) b(ld+k-j) +\cdots +
\sum_{j=(l-1)d}^{ld-1} a(j) b(ld+k-j) + \sum_{j=ld}^{ld+k} a(j) b(ld+k-j) .
\]
We now handle the individual terms in the foregoing identity. Let $r_j$ be the remainder of the
division of $j$ by $d$ for $j=0,1,\ldots, ld+k$. Then for $m=0,1,\ldots, l-1$ and $md\leq j \leq (m+1)d-1$,
we have $j = md +r_j$ with $0\leq r_j<d$. Then
\begin{equation}\label{Eq:lem-3}
\sum_{j=md}^{(m+1)d-1} a(j)b(ld+k-j)
\end{equation}
\[
=
\sum_{j: r_j\leq k} b(md+r_j) b\big( (l-m)d + k-r_j \big) 
+ \sum_{j: r_j > k} a(md+ r_j) b \big( (l-m)d + k-r_j \big) 
\]
\[
=
\sum_{j= 0}^{k} a(md+j) b\big( (l-m)d + k-j \big) 
 + \sum_{j: r_j >k} a(r_j) b\big( (l-m)d + k- r_j \big) 
\]
\[
=
a(md)b\big( (l-m)d \big) \sum_{j=0}^k a(j) b(k-j)  
 + a(md) b\big( (l-m)d \big) \sum_{j: r_j >k}
a(r_j) \fr{b\big( (l-m)d + k-r_j \big)}{b\big( (l-m)d \big)}
\]
We now claim that
\begin{equation}\label{Eq:lem-5}
\sum_{j: r_j >k}
a(r_j) \fr{b\big( (l-m)d + k-r_j \big)}{b\big( (l-m)d \big)}
=\sum_{r_j =k+1}^{d-1} a(r_j)\fr{b\big( (l-m)d + k-r_j \big)}{b\big( (l-m)d \big)} = 0.
\end{equation}
Note first that
the claim is clear if $k>\fr{d-1}{2}$ since $a(r_j)=0$ for $\fr{d-1}{2}<r_j\leq d-1$ by assumption.
By the same assumption, it easy to see that the terms in the foregoing sum for which
 $r_j> \fr{d-1}{2} \geq k$ vanish. Now suppose that $k<r_j \leq \fr{d-1}{2}$. It follows that
 $r_j-k < \fr{d-1}{2}$ and so, $d+k-r_j > \fr{d-1}{2}$. Thus we get
\[
a(r_j) \fr{b\big( (l-m)d + k-r_j \big)}{b\big( (l-m)d \big)}
=
a(r_j) \fr{b\big((l-m-1)d \big)}{b\big( (l-m)d \big)}
\fr{ b\big((l-m-1)d +d+k-j \big)}{b\big((l-m-1)d \big)} 
\]
\[
= a(r_j) \fr{b\big((l-m-1)d \big)}{b\big( (l-m)d \big)} b\big(d+k-j \big) 
= 0.
\]
 This proves the claim.
 Similarly, we have
\begin{equation}\label{Eq:lem-6}
\sum_{j=ld}^{ld+k} a(j) b(ld+k-j) = a(ld)\sum_{j=0}^k a(j)b(k-j).
\end{equation}
Now combine the relations (\ref{Eq:lem-2})--(\ref{Eq:lem-6}) to deduce that
\[
c(ld+k) = c(k) \sum_{m=0}^{l-1} a(md) b\big( (l-m)d \big) + c(k) a(ld)
= c(k)\sum_{m=0}^{l} a(md) b\big( (l-m)d \big)
\]
which gives the desired formula for the case $N=2$. Suppose now that the statement is true for $N-1$ and let
\[
s(k) = \sum_{i_2 =0}^{k} \sum_{i_3 =0}^{k-i_2}\cdots \sum_{i_{N-1}=0}^{k-i_2-\cdots-i_{N-1}}
z_{2}(i_2)\cdots z_{N-1}(i_{N-1}) z_N (k-i_2-\cdots-i_{N-1}).
\]
We have
\begin{equation}\label{triple}
 t(ld+k)
= \sum_{i_1=0}^{ld+k}\sum_{i_2=0}^{ld+k-i_1} \cdots \sum_{i_N=0}^{ld+k-i_1-\cdots - i_{N-1}}
z_1(i_1) \cdots z_{N-1}(i_{N-1}) z_N (ld+k-i_1-\cdots - i_{N-1}).
\end{equation}
\[
= \sum_{i_1=0}^{d-1}\sum_{i_2=0}^{ld+k-i_1} \cdots \sum_{i_N=0}^{ld+k-i_1-\cdots - i_{N-1}}
z_1(i_1) \cdots z_{N-1}(i_{N-1}) z_N (ld+k-i_1-\cdots - i_{N-1})
\]
\[
+ \sum_{i_1=d}^{2d-1}\sum_{i_2=0}^{ld+k-i_1} \cdots \sum_{i_N=0}^{ld+k-i_1-\cdots - i_{N-1}}
z_1(i_1) \cdots z_{N-1}(i_{N-1}) z_N (ld+k-i_1-\cdots - i_{N-1})
\]
\[
+\cdots + \sum_{i_1=(l-1)d}^{ld-1}\sum_{i_2=0}^{ld+k-i_1} \cdots \sum_{i_N=0}^{ld+k-i_1-\cdots - i_{N-1}}
z_1(i_1) \cdots z_{N-1}(i_{N-1}) z_N (ld+k-i_1-\cdots - i_{N-1})
\]
\[
+ \sum_{i_1=ld}^{ld+k}\sum_{i_2=0}^{ld+k-i_1} \cdots \sum_{i_N=0}^{ld+k-i_1-\cdots - i_{N-1}}
z_1(i_1) \cdots z_{N-1}(i_{N-1}) z_N (ld+k-i_1-\cdots - i_{N-1}).
\]
Let $r$ be the remainder of $i_1$ divided by $d$ for $i_1=0,1\ldots,ld+k$. Clearly for
$m=0,1,\ldots,l-1$ and $md\leq j\leq (m+1)d-1$ we have $i_1=md+r$ with $0\leq r <d$.
Then for $m=0,1,\ldots,l-1$ we have
\begin{equation}\label{term-m}
\sum_{i_1=md}^{(m+1)d-1}\sum_{i_2=0}^{ld+k-i_1} \cdots \sum_{i_N=0}^{ld+k-i_1-\cdots - i_{N-1}}
z_1(i_1) \cdots z_{N-1}(i_{N-1}) z_N (ld+k-i_1-\cdots - i_{N-1})
\end{equation}
\[
= \sum_{i_1=md}^{(m+1)d-1} z_1(md+r)\sum_{i_2=0}^{ld+k-i_1} \cdots \sum_{i_N=0}^{ld+k-i_1-\cdots - i_{N-1}}
z_2(i_2) \cdots z_{N-1}(i_{N-1})
\]
\[
\times z_N \big( (l-m)d+k-r- i_2-\cdots - i_{N-1}) \big)
\]
\[
= z_1(md)  \sum_{j=md}^{(m+1)d-1} z_1(r)\sum_{i_2=0}^{ld+k-i_1} \cdots \sum_{i_N=0}^{ld+k-i_1-\cdots - i_{N-1}}
z_2(i_2) \cdots z_{N-1}(i_{N-1})
\]
\[
\times z_N \big( (l-m)d+k-r- i_2-\cdots - i_{N-1}) \big)
\]
\[
= z_1(md)  \sum_{\substack{j=md\\ r\leq k}}^{(m+1)d-1} z_1(r)
\sum_{i_2=0}^{ld+k-i_1} \cdots \sum_{i_N=0}^{ld+k-i_1-\cdots - i_{N-1}}
z_2(i_2) \cdots z_{N-1}(i_{N-1})
\]
\[
\times z_N \big( (l-m)d+k-r- i_2-\cdots - i_{N-1}) \big)
\]
\[
 + z_1(md)  \sum_{\substack{j=md\\ r > k}}^{(m+1)d-1} z_1(r)
\sum_{i_2=0}^{ld+k-i_1} \cdots \sum_{i_N=0}^{ld+k-i_1-\cdots - i_{N-1}}
z_2(i_2) \cdots z_{N-1}(i_{N-1})
\]
\[
\times z_N \big( (l-m)d+k-r- i_2-\cdots - i_{N-1}) \big).
\]
Let $A$ and $B$ respectively  be the first and the second term on the right hand-side of the foregoing
identity. Then by the induction hypothesis
\[
A= z_1(md)  \sum_{i_1=0}^{k} z_1(i_1)
\sum_{i_2=0}^{(l-m)d+k-i_1} \cdots \sum_{i_N=0}^{(l-m)d+k-i_1-\cdots - i_{N-1}}
z_2(i_2) \cdots z_{N-1}(i_{N-1})
\]
\begin{equation}
\qquad \qquad \qquad \qquad \qquad \qquad \times z_N \big( (l-m)d+k-i_1- i_2-\cdots - i_{N-1}) \big)
\end{equation}
\[
= z_1(md) \sum_{i_1=0}^{k} z_1(i_1) s\big( (l-m)d + k-i_1 \big)
\]
\[
= z_1(md)  \sum_{i_1=0}^{k} z_1(i_1) s(k-i_1) \sum_{i_2=0}^{l-m} \sum_{i_3=0}^{l-m-i_2}
\cdots \sum_{i_N=0}^{(l-m)d+k-i_2-\cdots - i_{N-1}}
z_2(i_2 d) \cdots z_{N-1}(i_{N-1} d)
\]
\[
 \qquad \qquad \qquad \qquad \qquad \qquad \times z_N \big((l-m- i_2-\cdots - i_{N-1})d \big)
 \]
\[
= z_1(md) t(k) \sum_{i_2=0}^{l-m} \sum_{i_3=0}^{l-m-i_2}
\cdots \sum_{i_N=0}^{(l-m)d+k-i_2-\cdots - i_{N-1}} z_2(i_2 d) \cdots z_{N-1}(i_{N-1} d)
\]

As for $B$, we find by the induction hypothesis and~\eqref{Eq:key-1}
\[
\begin{split}
B&= z_1(md)  \sum_{\substack{i_1=md\\ r > k}}^{(m+1)d-1} z_1(r)
\sum_{i_2=0}^{ld+k-i_1} \cdots \sum_{i_N=0}^{ld+k-i_1-\cdots - i_{N-1}} z_2(i_2) \cdots z_{N-1}(i_{N-1}) \\
&  \qquad \qquad \qquad \qquad \qquad \qquad \times z_N \big( (l-m)d+k-r- i_2-\cdots - i_{N-1}) \big) \\
&= z_1(md)  \sum_{i_1=k+1}^{d-1} z_1(i_1) \sum_{i_2=0}^{(l-m)d+k-i_1} \cdots \sum_{i_N=0}^{(l-m)d+k-i_1-\cdots - i_{N-1}} z_2(i_2) \cdots z_{N-1}(i_{N-1}) \\
&\qquad \qquad \qquad \qquad \qquad \qquad \times z_N \big( (l-m)d+k-i_1- i_2-\cdots - i_{N-1}) \big) \\
&=  z_1(md)  \sum_{i_1=k+1}^{d-1} z_1(i_1) \sum_{i_2=0}^{(l-m-1)d+d+k-i_1} \cdots \sum_{i_N=0}^{(l-m-1)d+d+k-i_1-\cdots - i_{N-1}} z_2(i_2) \cdots z_{N-1}(i_{N-1}) \\
&\qquad \qquad \qquad \qquad \qquad \qquad \times z_N \big( (l-m-1)d+d+k-i_1- i_2-\cdots - i_{N-1}) \big) \\
&= z_1(md)  \sum_{i_1=k+1}^{d-1} z_1(i_1) s\big( (l-m-1)d + d+k-i \big) \\
&= z_1(md)  \sum_{i_1=k+1}^{\fr{d-1}{N}} z_1(i_1) \sum_{i_2=0}^{l-m-i_1} \cdots
\sum_{i_N=0}^{l-m-i_1-\cdots-i_{N-1}}  z_1(i_2) \cdots z_{N-1} (i_{N-1}) z_N(d+k-i_1-\cdots-i_{N-1}) \\
& \qquad \qquad \qquad \qquad \qquad \qquad \times \sum_{i_1=k+1}^{d-1} z_1(i_1) s(d+k-i_1) \\
&= z_1(md)  \sum_{i_1=k+1}^{\fr{d-1}{N}} z_1(i_1) \sum_{i_2=0}^{\fr{d-1}{N}} \cdots
 \sum_{i_N=0}^{\fr{d-1}{N}}  z_1(i_2) \cdots z_{N-1} (i_{N-1})
z_N(d+k-i_1-\cdots-i_{N-1}) \\
& \qquad \qquad \qquad \qquad \qquad \qquad \times \sum_{i_1=k+1}^{d-1} z_1(i_1) v(d+k-i_1) \\
\end{split}
\]
As $k-i_1, i_2,\ldots, i_{N-1} \leq (d-1)/N$, we deduce that $d-1\geq d+k-i_1-\cdots - i_{N-1} > (d-1)/N$
and therefore $z_N (d+k-i_1-\cdots - i_{N-1} ) =0$. It follows that
\[
s(d+k-i_1)
= \sum_{i_2=0}^{\fr{d-1}{n}} \cdots
\sum_{i_N=0}^{\fr{d-1}{N}}  z_1(i_2) \cdots z_{N-1} (i_{N-1}) z_N(d+k-i_1-\cdots-i_{N-1})
=0
\]
and hence $B=0$.
Regarding the last term of~\eqref{triple}, we similarly find
\[
\sum_{i_1=ld}^{ld+k}\sum_{i_2=0}^{ld+k-i_1} \cdots \sum_{i_N=0}^{ld+k-i_1-\cdots - i_{N-1}}
z_1(i_1) \cdots z_{N-1}(i_{N-1}) z_N (ld+k-i_1-\cdots - i_{N-1})
\]
\begin{equation}\label{l-term}
= z_1(ld) t(k).
\end{equation}
Then by~\eqref{triple}-\eqref{l-term},
we deduce after replacing $m$ with $i_1$ that
\[
t(ld+k) = t(k) \sum_{i_1=0}^{l} \sum_{i_2 =0}^{l-i_1} \cdots
\sum_{i_{N}=0}^{l-i_1-\cdots-i_{N-1}}
z_{1}(i_1 d)\cdots z_{N-1}(i_{N-1}d) z_N \big( (l-i_1-\cdots-i_{N-1})d \big).
\]
This completes the proof.

\section{Final comments}\label{sec:comments}
Note that in Theorem~\ref{thm:appl-5} while the sequence $\{\alpha_q(k)\}_{k=0}^\infty$ with the choices $x=q^{-n}$ or $x=q^n$
satisfies~\eqref{Eq:key-1} for $N=2$, it does not satisfy this condition for $N>2$. This justifies the fact that while
for any prime $p>1$ we have
\[
\sum_{k=0}^{p-1}\sum_{j=0}^k
\fr{{2j\choose j}^2{4j\choose 2j}
{2(k-j)\choose k-j}^2{4(k-j)\choose 2(k-j)}}{2^{8k} 3^{2k}} (8j+1) \big(8(k-j)+1 \big) \equiv p^2 \pmod{p^3},
\]
it is in general not true that
\[
\sum_{k=0}^{p-1} \sum_{j=0}^k \sum_{i=0}^{k-j}
\fr{ {2j\choose j}^2{4j\choose 2j} {2i\choose i}^2{4i\choose 2i}
{2(k-j-i)\choose k-j-i}^2{4(k-j-i)\choose 2(k-j-i)}}{2^{8k} 3^{2k}}
(8i+1)(8j+1) \big(8(k-j-i)+1 \big)
\]
\[
\equiv p^3 \equiv 0 \pmod{p^3}.
\]
Related to this we have the following conjecture which is based on computational evidence.
\begin{conjecture}
If $p>3$ is a prime number, then
\[
\sum_{k=0}^{p-1} \sum_{j=0}^k \sum_{i=0}^{k-j}
\fr{ {2j\choose j}^2{4j\choose 2j} {2i\choose i}^2{4i\choose 2i}
{2(k-j-i)\choose k-j-i}^2{4(k-j-i)\choose 2(k-j-i)}}{2^{8k} 3^{2k}}
(8i+1)(8j+1) \big(8(k-j-i)+1 \big)
\equiv 0 \pmod{p^2}.
\]
\end{conjecture}
\bigskip
%
\noindent{\bf Data Availability Statement.}
The authors confirm that the manuscript is self-contained and that the data supporting the findings of this study are available within the article.

\end{document}